\mag=\magstep1
\documentclass[reqno]{amsart}
\usepackage{graphicx}
\usepackage{amssymb}
\def\theorem#1{\bigbreak\noindent{\bf Theorem {#1}. }}

\def\lemma#1{\bigbreak\noindent{\bf Lemma {#1}. }}

\def\conjecture#1{\bigbreak\noindent{\bf Conjecture {#1}. }}

\def\bgns{\it}

\def\fin{\rm\bigbreak}

\def\proof{\noindent{\it Proof. }}
\def\proofof#1{\noindent{\it Proof of {#1}. }}

\def\trans{\raise.8em\hbox{\sevenit t}}

\newcommand\h{\eta}
\newcommand\tht{\theta}

\newcommand\la{\lambda}
\newcommand\m{\mu}

\newcommand\s{\sigma}
\newcommand\ta{\tau}

\newcommand\f{\phi}

\newcommand\ps{\psi}

\newcommand\G{\Gamma}

\newcommand\F{\Phi}

\newcommand\bN{\mathbb{N}}

\newcommand\bP{\mathbb{P}}
\newcommand\bQ{\mathbb{Q}}
\newcommand\bR{\mathbb{R}}

\newcommand\bZ{\mathbb{Z}}

\newcommand\eqqx[1]{\begin{equation}#1\end{equation}}            

\newcommand\tarrx[2]{\left\{\begin{array}{#1}#2\end{array}\right.}

\newcommand\arrx[2]{\begin{array}{#1}#2\end{array}}

\newcommand\parrx[2]{\left(\begin{array}{#1}#2\end{array}\right)}

\newcommand\paren[1]{\left(#1\right)}

\newcommand\gammas[2]{M_{#1}^{#2}}                                  
\newcommand\ppgr{\tilde Q}                                          

\newcommand\lrarrow{\longrightarrow}

\newcommand\vs{\vspace{1ex}}
\newcommand\vhs{\vspace{.5ex}}

\begin{document}

\title[The group generated by the gamma functions]
{\LARGE The group generated by the gamma functions $\Gamma(ax+1)$,\\ 
and its subgroup of the elements converging to constants}
\author{KAZUTO ASAI}
\dedicatory{Center for Mathematical Sciences, University of Aizu, \vs\\
 Aizu-Wakamatsu, Fukushima 965-8580, Japan \vs\\
{\rm e-mail: k-asai@u-aizu.ac.jp} \vs\\
Tel. 0242-37-2644 (Office), 0242-37-2752 (Fax)} 
\address{Center for Mathematical Sciences, University of Aizu, Aizu-Wakamatsu, Fukushima 965-8580, Japan}
\email{k-asai@u-aizu.ac.jp}

\maketitle

\hrule \vs\vs
Abstract \\ 
Let $G$ be the multiplicative group generated by the gamma functions $\Gamma(ax+1)$ $(a=1,2,\dots)$, and $H$ be the subgroup of all elements of $G$ that converge to nonzero constants as $x\rightarrow\infty$. The quotient group $G/H$ is the group of equivalence classes of $G$, where $f$ and $g$ are equivalent $\iff f\sim Cg$ $(x\rightarrow\infty)$ for some $C\not=0$. We show that $G/H\simeq\mathbb{Q}^+$. A similar consideration is possible for the case that the gamma functions $\Gamma(ax+1)$ with $a\in\mathbb{R}^+$ are concerned, and we show that $G/H\simeq\mathbb{Z}\times\mathbb{R}\times\mathbb{R}$. 

Also, several concrete examples of the elements of $H$ are constructed, e.g., it holds that $\frac{\binom{18n}{12n,3n,3n}}{\binom{18n}{9n,8n,n}}\longrightarrow\sqrt{\frac{2}{3}}$ $(n\rightarrow\infty)$, where $\binom{*}{*,\dots,*}$ denotes a multinomial coefficient. 
\vskip2ex
\hrule
\bigskip

\section{Introduction}

Throughout this paper, let $\bP$ denote the set of all positive integers, $\bN$ denote the set of all nonnegative integers, $\bQ^+$ denote the multiplicative group of all positive rational numbers, and $\bR^+$ denote the multiplicative group of all positive real numbers. 

So far, we have found many kinds of formulas for the asymptotic behavior of the gamma function \cite{gosp:hyper78, liu:newv07, mortici:ulti-fact09, mortici:subs11, mortici:monot12}, including Ramanujan's factorial approximation and its analogues \cite{mortici:raman11, mortici:ramest11, mortici:raman12, raman:lost88}. The origin of them is Stirling's formula discovered by A. de Moivre and J. Stirling \cite{stirl:methd30, twed-stir:methd03}:
\eqqx{\G(x+1)\sim\sqrt{2\pi x}\paren{\frac{x}{e}}^x\qquad(x\lrarrow\infty).\label{stir-approx}}
In this paper, we use this simple formula to evaluate the asymptotic behavior of the ratio of the products of the gamma functions. 

Let $G$ be the multiplicative group generated by the gamma functions $\G(ax+1)$ $(a\in\bP)$, and $H$ be the subgroup of all elements of $G$ that converge to nonzero constants as $x\rightarrow\infty$. By using the notation: 
\eqqx{\gammas{b_1,\dots,b_t}{a_1,\dots,a_s}=\gammas{b_1,\dots,b_t}{a_1,\dots,a_s}(x)
=\frac{\prod_{k=1}^s\G(a_kx+1)}{\prod_{k=1}^t\G(b_kx+1)},}
we have
\eqqx{G=\{\gammas{b_1,\dots,b_t}{a_1,\dots,a_s}\mid a_1,\dots,a_s,
b_1,\dots,b_t\in\bP;\,s,t\in\bN\}.}
Here, when $st=0$, $\gammas{b_1,\dots,b_t}{a_1,\dots,a_s}$ becomes $\gammas{}{a_1,\dots,a_s}$, $\gammas{b_1,\dots,b_t}{}$, or $\gammas{}{}=1$. 

We consider the quotient group $G/H$, which is the group of equivalence classes of $G$, where $f$ and $g$ are defined to be equivalent when $f\sim Cg$ $(x\rightarrow\infty)$ for some nonzero constant $C$. We show the following. 

\theorem1\bgns 
It holds that \ $G\simeq G/H\simeq\bQ^+$.
\fin

A similar consideration is possible for the case that the gamma functions $\G(ax+1)$ with $a\in\bR^+$ are taken as the generators. Let $\tilde G$ be the multiplicative group generated by $\G(ax+1)$ $(a\in\bR^+)$, and $\tilde H$ be the subgroup of all elements of $\tilde G$ that converge to nonzero constants as $x\rightarrow\infty$. 

\theorem2\bgns
It holds that \ $\tilde G/\tilde H\simeq\bZ\times\bR\times\bR$. 
\fin

In Section~3, we study concrete elements of $H$ in a combinatorial context. For partitions $\la,\m$, a primitive solution $(\la;\m)$ to the condition for $M_{\m}^{\la}\in H$ is defined, and we prove that there exists a primitive solution of length exceeding $n$ for every integer $n$. Several examples of the primitive solutions are also given. 

\bigskip

\section{Homomorphisms and proofs}

First we give the asymptotic representation of the elements of $G$: 

\lemma1\bgns
\eqqx{\arrx{l}{\gammas{b_1,\dots,b_t}{a_1,\dots,a_s}(x)\sim\sqrt{\frac{a_1\dots a_s}
{b_1\dots b_t}}\,(2\pi x)^{\frac{s-t}{2}}\paren{\frac{x}{e}}^{
x\paren{\sum_{k=1}^sa_k-\sum_{k=1}^tb_k}}
\paren{\frac{a_1^{a_1}\dots a_s^{a_s}}{b_1^{b_1}\dots b_t^{b_t}}}^x. \vs\\
(x\lrarrow\infty)}
\label{gammas-limit}}\fin

\proof This is a direct consequence of Stirling's formula \eqref{stir-approx}. Put $x\rightarrow ax$ and calculate $\gammas{b_1,\dots,b_t}{a_1,\dots,a_s}(x)$. 
\qed\fin

By this lemma, we have the condition for the elements of $G$ to be contained in $H$. Indeed, $\gammas{b_1,\dots,b_t}{a_1,\dots,a_s}(x)$ converges to a nonzero constant if and only if the three factors $(2\pi x)^{\frac{s-t}{2}}$, $\paren{\frac{x}{e}}^{x*}$, $(*{*}*)^x$ of \eqref{gammas-limit} are constants equal to 1. 

\lemma2\bgns For $\gammas{b_1,\dots,b_t}{a_1,\dots,a_s}\in H$, it is necessary and sufficient that (i) $s=t$, (ii) $\sum_{k=1}^sa_k=\sum_{k=1}^tb_k$, and (iii) $a_1^{a_1}\dots a_s^{a_s}=b_1^{b_1}\dots b_t^{b_t}$. 
\fin

Hereafter, we consider homomorphisms from $G$ to certain groups for preparation of the proof of Theorem~1. Let $\ppgr$ denote the multiplicative group generated by $\{p^p\mid p:\mbox{a prime}\}$. Let $\f_i$ $(i=1,2,3)$ be homomorphisms defined below: 
\eqqx{\arrx{l}{
\f_1:G\lrarrow\bZ:\gammas{b_1,\dots,b_t}{a_1,\dots,a_s}\longmapsto s-t \\
\f_2:G\lrarrow\bZ:\gammas{b_1,\dots,b_t}{a_1,\dots,a_s}\longmapsto 
\sum_{k=1}^sa_k-\sum_{k=1}^tb_k \\
\f_3:G\lrarrow\ppgr:\gammas{b_1,\dots,b_t}{a_1,\dots,a_s}\longmapsto 
\frac{a_1^{a_1}\dots a_s^{a_s}}{b_1^{b_1}\dots b_t^{b_t}}.}
\label{homomorphisms}}
This definition is possible because each element of $G$ is uniquely expressed by the symbol $\gammas{b_1,\dots,b_t}{a_1,\dots,a_s}$ up to a permutation of indices and a cancellation of identical upper/lower indices. Then we have a homomorphism: 
\eqqx{\F:G\lrarrow\bZ\times\bZ\times\ppgr:g\longmapsto(\f_1(g),\f_2(g),\f_3(g)).}

\lemma3\bgns $\F$ is a surjection. 
\fin

\proof To begin with we note that $\f_3$ is well defined, say, $\f_3(G)\subset\ppgr$. Take any $g=\gammas{b_1,\dots,b_t}{a_1,\dots,a_s}\in G$. If $a_1=p^em$ and $p\nmid m$ for a prime factor $p$, we have $a_1^{a_1}=(p^em)^{p^em}=p^{p(p^{e-1}em)}m^{a_1}$. Hence there exist prime numbers $p_1,\dots,p_l$, and integers $e_1,\dots,e_l$ such that
\eqqx{\frac{a_1^{a_1}\dots a_s^{a_s}}{b_1^{b_1}\dots b_t^{b_t}}=
\prod_{k=1}^lp_k^{p_ke_k},
\label{phi-3}}
and therefore $\f_3(G)\subset\ppgr$. Next we confirm that  $\f_3$ is a surjection. For an arbitrary element $y=\frac{p_1^{p_1}\dots p_s^{p_s}}{q_1^{q_1}\dots q_t^{q_t}}$ of $\ppgr$ with prime numbers $p_1,\dots,p_s,q_1,\dots,q_t$ (repetition allowed), taking $g=\gammas{q_1,\dots,q_t}{p_1,\dots,p_s}$, we have $y=\f_3(g)$. 

Now, we prove $\F$ is a surjection. Let $(d,l,y)$ be an arbitrary element of $\bZ\times\bZ\times\ppgr$. Since $\f_3$ is a surjection, we have $y=\f_3(g)$ for some $g=\gammas{b_1,\dots,b_t}{a_1,\dots,a_s}\in G$. If $\f_1(g)\not=d$, we can take $g_1=\gammas{b_1,\dots,b_t}{a_1,\dots,a_s,1,\dots,1}$ or $\gammas{b_1,\dots,b_t,1,\dots,1}{a_1,\dots,a_s}$ such that $\f_1(g_1)=d$, $\f_3(g_1)=y$. If $\f_2(g_1)=m\not=l$, consider $g_2=\gammas{6,2,1}{4,3,3}$. We see $\F(g_2)=(0,1,1)$. Thus, letting $g_3=g_1g_2^{l-m}$, we have $\F(g_3)=(d,l,y)$.  \qed
\fin

\proofof{Theorem~1}
(i) $G\simeq\bQ^+$: The mapping: $\ps:G\lrarrow\bQ^+$ defined by 
\eqqx{\ps(\gammas{b_1,\dots,b_t}{a_1,\dots,a_s})=
\frac{p_{a_1}\dots p_{a_s}}{p_{b_1}\dots p_{b_t}}}
is confirmed to be an isomorphism, where $p_i$ denotes the $i$-th prime. 

(ii) $G/H\simeq\bQ^+$: Apply the fundamental homomorphism theorem: $G/\ker\F\simeq\F(G)$ to the above-defined $\F:G\lrarrow\bZ\times\bZ\times\ppgr$. By definition of $\F$ and Lemma~2, $\ker\F=H$. Hence, together with Lemma~3, we have $G/H\simeq\bZ\times\bZ\times\ppgr$. 

There exist isomorphisms $i:\ppgr\lrarrow\bQ^+$ and $j:\bZ\times\bQ^+\lrarrow\bQ^+$ defined by $i(p^p)=p$ for every prime $p$ and
\eqqx{j(n,2^{n_2}3^{n_3}5^{n_5}\dots)=2^n3^{n_2}5^{n_3}\dots.}
Therefore 
\eqqx{\bZ\times\bZ\times\ppgr\simeq\bZ\times\bZ\times\bQ^+\simeq\bQ^+.}
This proves $G/H\simeq\bQ^+$.  \qed
\fin

Next, we extend the homomorphism $\F$ to $\tilde\F$ defined on $\tilde G$ in order to prove Theorem~2. For that purpose we extend the homomorphisms $\f_i$ to $\tilde\f_i$ defined on $\tilde G$ by the correspondence formulas used in \eqref{homomorphisms}. Then $\tilde \F$ is extended to a homomorphism: 
\eqqx{\tilde\F:\tilde G\lrarrow\bZ\times\bR\times\bR^+:g\longmapsto(
\tilde\f_1(g),\tilde\f_2(g),\tilde\f_3(g)).}

\lemma4\bgns $\tilde\F$ is a surjection. 
\fin

\proof Let $(d,x,y)$ be an arbitrary element of $\bZ\times\bR\times\bR^+$. Clearly $\tilde\f_3(g)=y$ for some $g\in\tilde G$.  In a similar manner as in the proof of Lemma~3, we have $\tilde\f_1(g_1)=d$ and $\tilde\f_3(g_1)=y$ for some $g_1\in\tilde G$. Now take $\h\in(1/e^{1/e},1)$, then $x^x=\h$ has distinct two real solutions $\tht_1,\tht_2\in(0,1)$. Let $\gammas{\tht_1}{\tht_2}=g(\h)$, then $\tilde\F(g(\h))=(0,\tht_1-\tht_2,1)$, where $\tht_1-\tht_2$ takes an arbitrary nonzero value in $(-1,1)$ depending on $\h$. Hence choosing suitable $\h_1,\dots,\h_l$, we have $g_2=g_1g(\h_1)\dots g(\h_l)$ such that $\tilde\F(g_2)=(d,x,y)$.  \qed
\fin

\proofof{Theorem~2}
Apply again the fundamental homomorphism theorem to $\tilde\F$. We have
\eqqx{\tilde G/\tilde H=\tilde G/\ker\tilde\F\simeq\tilde\F(\tilde G)=
\bZ\times\bR\times\bR^+\simeq\bZ\times\bR\times\bR.}
\qed\fin

\bigskip

\section{Partitions and primitive solutions}

In this section, we construct concrete examples of the elements of $H$, and study primitive solutions defined below corresponding to such elements $M_{b_1,\dots,b_s}^{a_1,\dots,a_s}(x)$ of $H$ that generate $H$ by ordinary multiplication and the variable transformation $x\lrarrow kx$ for every positive integer $k$. 

A partition of a positive integer $n$ is a weakly decreasing sequence of positive integers: $\la=(\la_1,\la_2,\dots,\la_s)$ such that $|\la|=\la_1+\la_2+\cdots+\la_s=n$. (See \cite{andr:parti76, andr:parti04} etc.) Each $\la_i$ is called a part of $\la$ and the integer $s$ is called the length of $\la$ denoted by $l(\la)$. Let $\la=(\la_1,\dots,\la_s)$ and $\m=(\m_1,\dots,\m_t)$ be partitions of length $s$ and $t$, respectively. Set $\gammas{\m}{\la}=\gammas{\m_1,\dots,\m_t}{\la_1,\dots,\la_s}$, then any element of $G$ is expressed in this form. Denote $\la^{\la}=\la_1^{\la_1}\dots\la_s^{\la_s}$; $k\la=(k\la_1,\dots,k\la_s)$ for a positive integer $k$; and denote by $\la\oplus\tilde\la$, the rearrangement of $(\la_1,\dots,\la_s,\tilde\la_1,\dots,\tilde\la_{s'})$ in decreasing order. The condition for $\gammas{\m}{\la}\in H$ in Lemma~2 is rewritten as
\eqqx{\mbox{(i)}\ l(\la)=l(\m)\qquad\mbox{(ii)}\ |\la|=|\m|\qquad
\mbox{(iii)}\ \la^{\la}=\m^{\m}.\label{parti-eq}}
Equations \eqref{parti-eq} have always solutions $\la=\m$, which we call trivial solutions. For every solution $(\la;\m)$ to \eqref{parti-eq}, $l(\la)=l(\m)$ and $|\la|=|\m|$ are called the length and the size of the solution, respectively. The solutions $(\la;\m)$ and $(\m;\la)$ are usually identified. If $(\la;\m)$ is a solution, then $k(\la;\m)=(k\la;k\m)$ is also a solution for every positive integer $k$, because 
\eqqx{(k\la)^{k\la}=k^{k|\la|}\paren{\la^{\la}}^k=k^{k|\m|}\paren{\m^{\m}}^k
=(k\m)^{k\m}.}
These solutions are called equivalent to each other. In addition, if two solutions $(\la;\m)$, $(\tilde\la;\tilde\m)$ of positive lengths exist, then $(\la;\m)\oplus(\tilde\la;\tilde\m)=(\la\oplus\tilde\la;\m\oplus\tilde\m)$ is also a solution, decomposable into two solutions. Hence it is important to find nontrivial solutions that can not be written in the form $k(\la;\m)$ $(k\ge2)$ nor $(\la;\m)\oplus(\tilde\la;\tilde\m)$, which we call primitive solutions. It is easily seen that there are no primitive solutions of length $\le2$. (For $n=1$, clear. For $n=2$, the function $f(x)=(a+x)^{a+x}(a-x)^{a-x}$ $(x\in[0,a])$ is strictly increasing, and so there is no solution except trivial $(\la;\la)$.)

\theorem3\bgns
For every positive integer $n$, there exists a primitive solution to \eqref{parti-eq} of length exceeding $n$. 
\fin
\proof For convenience, we sometimes use the notation $\binom{\la}{\m}$ for a solution $(\la;\m)$ to \eqref{parti-eq}. Also, we denote by $\{\la\}$ the multiset which consists of all parts of $\la$. We prove that the following is a primitive solution to \eqref{parti-eq} for $n\ge8$:
\eqqx{\parrx{lllll}{2^n,&2^{n-2},2^{n-2},&\overbrace{2,2,\dots,2}^{2^{n-2}},
&\overbrace{2,2,\dots,2}^{2^{n-2}},&\overbrace{2,2,\dots,2}^{2^{n-2}} \\
2^{n-1},&2^{n-1},2^{n-1},&\underbrace{4,4,\dots,4}_{2^{n-2}},
&\underbrace{1,1,\dots,1}_{2^{n-2}},&\underbrace{1,1,\dots,1}_{2^{n-2}}}.
\label{p-solution}}
One can confirm that \eqref{p-solution} is a solution of length $3\times2^{n-2}+3$ and of size $3\times2^n$, and that $\la^{\la}=\m^{\m}=2^{2^{n-1}(3n+1)}$. It is also clear that \eqref{p-solution} is not a $k$ times multiple of some solution for $k\ge2$. Thus it suffices to show that \eqref{p-solution} is not decomposable into two solutions. 

Suppose \eqref{p-solution} is decomposed into $(\s;\ta)\oplus(\tilde\s;\tilde\ta)$. Write \eqref{p-solution} as $\binom{\la}{\m}=\binom{\la^1,\la^2}{\m^1,\m^2}$, where $\la^1$ is the first three parts of $\la$, $\la^2$ is the rest of it, and $\m^1,\m^2$ are defined similarly. If $\s$ and $\ta$ are composed by choosing only the parts of $\la^2$ and $\m^2$, respectively, then from $|\s|=|\ta|$, it follows that $\binom{\s}{\ta}$ consists of the blocks $\binom{2,2,2}{4,1,1}$, which contradicts $\s^{\s}=\ta^{\ta}$. For the case that $\s$ and $\ta$ contain only some parts of $\la^1$ and $\m^2$ (or $\la^2$ and $\m^1$), respectively, the only possibility is $n=2,4$ ($n=2$). Also, it is clearly impossible that $\s$ and $\ta$ could contain only some parts of $\la^1$ and $\m^1$, respectively. Thus we should deal with the case that $\s$ contains both parts of $\la^1$ and $\la^2$ or $\ta$ contains both parts of $\m^1$ and $\m^2$. If $\s$ contains all parts of $\la^1$ and $\ta$ contains no parts of $\m^1$, we have
\eqqx{\frac{\s^{\s}}{\ta^{\ta}}\ge\frac{2^{n2^{n}}2^{(n-2)2^{n-2}\times2}}{2^{8\times2^{n-2}}}
=2^{(3n-6)2^{n-1}}>1}
for $n\ge3$. The alternative case that $\s$ contains no parts of $\la^1$ and $\ta$ contains all parts of $\m^1$ is very similar. 
Hence we consider the case that $\s$ or $\ta$ has a nonempty proper submultiset of $\{\la^1\}$ or $\{\m^1\}$ as parts, respectively. (The other cases already appear above for $(\s;\ta)$ or $(\tilde\s;\tilde\ta)$.) 
If $\{\s\}\cap\{\la^1\}=\{2^{n-2}\}$ and $\{\ta\}\cap\{\m^1\}=\emptyset$ (as multiset),  then 
\eqqx{\frac{\s^{\s}}{\ta^{\ta}}\ge\frac{2^{(n-2)2^{n-2}}2^{2\times(2^{n-2}-1)}}{2^{8\times2^{n-2}}}
=2^{(n-8)2^{n-2}-2}.}
Hence for $n\ge9$, $\frac{\s^{\s}}{\ta^{\ta}}>1$, and so $(\s;\ta)$ is not a solution. For $n=8$, the only possibility that fits $\frac{\s^{\s}}{\ta^{\ta}}=1$ is
\eqqx{(2^{n-2},\underbrace{2,2,\dots,2}_{2^{n-2}}\,;\,
\underbrace{4,4,\dots,4}_{2^{n-2}},1),
\label{e1-solution}}
but $|\s|\not=|\ta|$. Therefore, for $n\ge8$, $(\s;\ta)$ is not a solution. 

If $\{\s\}\cap\{\la^1\}=\{2^{n-2}\}$ and $\{\ta\}\cap\{\m^1\}=\{2^{n-1}\}$,  then 
\eqqx{\frac{\s^{\s}}{\ta^{\ta}}\le\frac{2^{(n-2)2^{n-2}}2^{2\times2^{n-2}\times3}}{2^{(n-1)2^{n-1}}}
=2^{(6-n)2^{n-2}}.}
Hence for $n\ge7$, $(\s;\ta)$ is not a solution. 

Although we can proceed in a similar manner, we give one more case: $\{\s\}\cap\{\la^1\}=\{2^{n-2},2^{n-2}\}$ and $\{\ta\}\cap\{\m^1\}=\{2^{n-1}\}$, that has a little different flavor. Let $\s^2$ be the partition consists of the parts of $\s$ contained in $\la^2$, and $\ta^2$ be defined similarly. Since $2^{n-2}+2^{n-2}=2^{n-1}$, we have $|\s^2|=|\ta^2|$. As $\frac{2^{(n-2)2^{n-2}}2^{(n-2)2^{n-2}}}{2^{(n-1)2^{n-1}}}=2^{-2^{n-1}}<1$, $\ta$ contains 1 parts, and by parity, at least two 1 parts. For the partition $\tilde\s^2$ obtained by exclusion of a 2 part from $\s^2$, and the partition $\tilde\ta^2$ obtained by exclusion of two 1 parts from $\ta^2$, we have $|\tilde\s^2|=|\tilde\ta^2|$ and $l(\tilde\s^2)=l(\tilde\ta^2)$. Hence $\binom{\tilde\s^2}{\tilde\ta^2}$ is composed of $s$ blocks of $\binom{2,2,2}{4,1,1}$ $(s\ge0)$, and therefore
\eqqx{\frac{\s^{\s}}{\ta^{\ta}}=2^{-2^{n-1}}\frac{2^22^{2\times3\times s}}{2^{8s}}=2^{-2^{n-1}-2(s-1)}=1.}
This fails for all $n\ge3$. 
\qed
\fin

The solution \eqref{p-solution} is primitive also for $n=6,7$, which is confirmed by showing the equations:
\eqqx{\tarrx{l}{a+b+c=d+e+f \\
2^na+2^{n-2}b+2c=2^{n-1}d+4e+f \\
2^{n2^na}2^{(n-2)2^{n-2}b}2^{2c}=2^{(n-1)2^{n-1}d}2^{8e}}}
have no solution $(a,b,c,d,e,f)$ of nonnegative integers with $a\le1$, $b\le2$, $c\le3\times2^{n-2}$, $d\le3$, $e\le2^{n-2}$, $f\le2^{n-1}$ for $n=6,7$, except $0$ or $(1,2,3\times2^{n-2},3,2^{n-2},2^{n-1})$. However, for $n=5$, \eqref{p-solution} is decomposed into:
\eqqx{\parrx{c}{8,2,2,2,2,2,2,2,2 \\ 4,4,4,4,4,1,1,1,1}\oplus
\parrx{c}{32,8,2,2,2,2,2,2,2,2,2,2,2,2,2,2,2,2 \\
16,16,16,4,4,4,1,1,1,1,1,1,1,1,1,1,1,1}.
\label{n=5-solution}}

By computational calculation, many primitive solutions of length $\ge4$ are easily found. (Several examples are listed below.) However, no solutions of length 3 are found except ones equivalent to $(12,3,3;9,8,1)$ for the size $\le2000$. 

\begin{table}[ht]
\begin{tabular}{|c|c|} \hline
length&primitive solutions \\ \hline
3&$\binom{12,3,3}{9,8,1}$ \\ \hline
4&$\binom{9,4,4,2}{6,6,6,1}$ \\ \hline
5&$\binom{10,4,2,2,2}{8,5,5,1,1}$, $\binom{8,3,3,3,3}{6,6,4,2,2}$,  
$\binom{16,6,3,3,2}{12,8,8,1,1}$, $\binom{14,6,4,3,3}{12,7,7,2,2}$, 
$\binom{12,5,5,4,4}{10,8,6,3,3}$ \\ \hline
6&$\arrx{c}{\binom{6,2,2,2,2,2}{4,4,3,3,1,1}, \binom{8,3,3,2,2,2}{6,4,4,4,1,1}, 
\binom{10,3,3,3,3,2}{6,6,5,5,1,1}, \binom{12,5,5,2,2,2}{10,6,6,4,1,1}, 
\binom{12,4,4,3,3,2}{8,6,6,6,1,1}, \vhs\\
\binom{10,6,3,3,3,3}{9,5,5,4,4,1}, \binom{10,4,4,4,3,3}{8,6,5,5,2,2}, 
\binom{10,9,4,2,2,2}{12,5,5,3,3,1}}$ \\ \hline
7&$\arrx{c}{\binom{9,2,2,2,2,2,2}{6,6,3,3,1,1,1}, \binom{12,2,2,2,2,2,2}{8,6,6,1,1,1,1}, 
\binom{9,4,4,4,2,2,2}{8,6,3,3,3,3,1}, \binom{15,3,2,2,2,2,2}{10,9,5,1,1,1,1}, 
\binom{12,6,2,2,2,2,2}{9,8,4,4,1,1,1}, \vhs\\
\binom{12,4,4,4,2,2,2}{8,8,6,3,3,1,1}, \binom{12,3,3,3,3,3,3}{9,6,6,4,2,2,1}, 
\binom{9,4,4,4,4,4,1}{8,6,6,3,3,2,2}}$ \\ \hline
8&$\binom{10,3,3,3,3,3,3,1}{9,5,5,2,2,2,2,2}$, $\binom{12,3,3,3,3,2,2,2}{9,6,4,4,4,1,1,1}$ 
\\ \hline
9&$\binom{8,2,2,2,2,2,2,2,2}{4,4,4,4,4,1,1,1,1}$, $\binom{10,3,3,2,2,2,2,2,2}{6,5,5,4,4,1,1,1,1}$ 
\\ \hline
10&$\binom{9,4,2,2,2,2,2,2,2,2}{8,3,3,3,3,3,3,1,1,1}$ \\ \hline
\end{tabular}
\bigskip
\caption{Primitive solutions of length $\le10$ and size $\le30$}\label{tab:prim-sol}
\end{table}

\conjecture1\bgns A primitive solution to \eqref{parti-eq} of length 3 is unique. 
\fin

\pagebreak

\bibliography{../rrefs/refgamma}

\begin{thebibliography}{10}

\bibitem{andr:parti76}
G.~E. Andrews.
\newblock {\em The Theory of Partitions}.
\newblock Cambridge University Press, Cambridge, 1976/1998.

\bibitem{andr:parti04}
G.~E. Andrews;~K. Eriksson.
\newblock {\em Integer Partitions}.
\newblock Cambridge University Press, Cambridge, 2004.

\bibitem{gosp:hyper78}
R.~W. Gosper.
\newblock Decision procedure for indefinite hypergeometric summation.
\newblock {\em Proc. Natl. Acad. Sci. USA}, 75(1):40--42, 1978.

\bibitem{liu:newv07}
Z.~Liu.
\newblock A new version of the {Stirling} formula.
\newblock {\em Tamsui Oxf. J. Math. Sci.}, 23(4):389--392, 2007.

\bibitem{mortici:ulti-fact09}
C.~Mortici.
\newblock An ultimate extremely accurate formula for approximation of the
  factorial function.
\newblock {\em Arch. Math. (Basel)}, 93(1):37--45, 2009.

\bibitem{mortici:raman11}
C.~Mortici.
\newblock On {Ramanujan}'s large argument formula for the gamma function.
\newblock {\em Ramanujan J.}, 26(2):185--192, 2011.

\bibitem{mortici:ramest11}
C.~Mortici.
\newblock {Ramanujan}'s estimate for the gamma function via monotonicity
  arguments.
\newblock {\em Ramanujan J.}, 25(2):149--154, 2011.

\bibitem{mortici:subs11}
C.~Mortici.
\newblock A substantial improvement of the {Stirling} formula.
\newblock {\em Appl. Math. Lett.}, 24(8):1351--1354, 2011.

\bibitem{mortici:monot12}
C.~Mortici.
\newblock Gamma function estimates via completely monotonicity arguments.
\newblock {\em Carpathian J. Math.}, 28(1):93--102, 2012.

\bibitem{mortici:raman12}
C.~Mortici.
\newblock An improvement of the {Ramanujan} formula for approximation of the
  {Euler} gamma function.
\newblock {\em Carpathian J. Math.}, 28(2):301--304, 2012.

\bibitem{raman:lost88}
S.~Ramanujan.
\newblock {\em The lost notebook and other unpublished papers. With an
  introduction by George E. Andrews}.
\newblock Springer-Verlag, Berlin; Narosa Publishing House, New Delhi, 1988.

\bibitem{stirl:methd30}
J.~Stirling.
\newblock {\em Methodus differentialis, sive tractatus de summation et
  interpolation serierum infinitarium}.
\newblock London, 1730.

\bibitem{twed-stir:methd03}
I.~Tweddle.
\newblock {\em {James Stirling}'s methodus differentialis: An annotated
  translation of {Stirling}'s text (Sources and studies in the history of
  mathematics and physical sciences)}.
\newblock Springer-Verlag, London, 2003.

\end{thebibliography}
\bibliographystyle{plain}


\end{document}